\documentclass{amsart}
\usepackage{amsmath,amsthm,amssymb} 
\usepackage{color, a4}
\usepackage{graphicx}

\begin{document}
\newtheorem{theorem}{Theorem}[section]
\newtheorem{definition}[theorem]{Definition}
\newtheorem{corollary}[theorem]{Corollary}
\newtheorem{remark}[theorem]{Remark}
\newtheorem{example}[theorem]{Example}
\newtheorem{lemma}[theorem]{Lemma}
\newtheorem{proposition}[theorem]{Proposition}

\newcommand{\R}{{\mathbb R}}
\newcommand{\C}{{\mathbb C}}
\newcommand{\Z}{{\mathbb Z}}
\newcommand{\In}{\mbox{\rm Inn}}
\newcommand{\Aut}{\mbox{\rm Aut}}
\newcommand{\id}{\mbox{\rm id}}

\title{Quandles and symmetric quandles for higher dimensional knots}

\author{Seiichi Kamada}
\address{Department of Mathematics, Osaka City University, 
Osaka 558-8585, Japan} 
\email{skamada@sci.osaka-cu.ac.jp} 

\thanks{The author was partially supported by KAKENHI (21340015, 23654027). }

\date{}
\maketitle

\begin{abstract}
A symmetric quandle is a quandle with a good involution.  For a knot in $\R^3$, a knotted surface in $\R^4$ or an $n$-manifold knot in $\R^{n+2}$, the knot symmetric quandle is defined.  We introduce the notion of a symmetric quandle presentation, and show how to get a presentation of a knot symmetric quandle from a diagram. 
\end{abstract}

\section{Introduction}

A {\it rack} is a set $X$ with a binary operation $X \times X \to X; 
(a,b) \mapsto a^b$ satisfying that 
(1) for any $x$ and $y$ of $X$, there is a unique element $z$ with 
$z^y= x$, and 
(2) for all $x,y,z \in X$, $(x^y)^z = (x^z) ^{(y^z)}$. 
A {\it quandle} is a rack satisfying 
that (3) $x^x=x$ for all $x \in X$. 
(The condition (2) is called (right) self-distributivity.  Refer to \cite{Dehornoy2000} for details on algebras with this property.)
In \cite{Joyce1982, Matveev1982}, the {\it knot quandle} (or the {\it fundamental quandle}) of a knot is defined and it is proved that the knot quandle is a complete invariant of a classical knot up to orientations of knots and the ambient space.  

These notions are defined not only for classical knots but also for (PL and locally flat, or smooth) proper and oriented $n$-submanifolds of an oriented $(n+2)$-manifold (cf. \cite{FennRourke1992, FRS2003}). 
Quandles are very useful for studying oriented knots in $\R^3$ and oriented surface knots in $\R^4$ (cf.\cite{CJKLS03}).  In \cite{FennRourke1992} the notion of a presentation of a rack or a quandle is established and it is shown that a presentation of the knot quandle of an oriented knot is obtained from a diagram (cf. \cite{Joyce1982, Matveev1982}).  
In this paper we show how to obtain a presentation of the knot quandle of a closed oriented $n$-submanifold of $\R^{n+2}$. 

In \cite{Kamada2007, KamadaOshiro2010} the notions of symmetric quandles, knot symmetric quandles (or the fundamental symmetric quandles), 
colorings and quandle homological invariants using symmetric quandles 
for knots in $\R^3$ and surface knots in $\R^4$ were introduced and studied.  However presentations of symmetric quandles were not explicitly discussed.  
The purpose of this paper is to establish the notion of presentations of symmetric quandles, as an analogy of that of presentations of quandles due to \cite{FennRourke1992}.  We define the notion of the  {\it knot symmetric quandle} 
of a  proper $n$-submanifold of an $(n+2)$-manifold, and 
show how to get a presentation of the knot symmetric quandle of a closed $n$-submanifold of $\R^{n+2}$ from its diagram.

\section{Symmetric quandles} 

Let $X$ be a quandle or a rack.  We follow the notation due to Fenn and Rourke \cite{FennRourke1992}. For example, $x^{y^{-1}}$ means the unique element $z$ with $z^y=x$.  

\begin{definition}[\cite{Kamada2007, KamadaOshiro2010}]{\rm 
An involution $\rho:X \to X$ is  a {\it good involution} if 
$$
\rho(x^y) = \rho(x)^y \quad \mbox{and} \quad x^{\rho(y)} = x^{y^{-1}} 
$$
for every $x$ and $y$ of $X$.  
A {\it symmetric quandle} (or a {\it symmetric rack}, resp.) is a pair $(X, \rho)$ of a quandle (or a rack, resp.) $X$ and a good involution $\rho: X \to X$. 
}\end{definition}

\begin{example}{\rm 
Let $X$ be a trivial quandle, i.e., $x^y = x$ for all $x, y \in X$.  Then every involution of $X$ is a good involution. 
Conversely, a quandle such that every involution is a good involution is a trivial quandle (cf. \cite{KamadaOshiro2010}). 
}\end{example}

\begin{example}{\rm 
Let $X$ be a kei (\cite{Takasaki1943}), i.e, it is a quandle satisfying $(x^y)^y = x$ for all $x, y \in X$. 
Then the identity map $\rho: X \to X$  is a good involution.   
Conversely a quandle such that the identity map is a good involution is a kei (cf. \cite{KamadaOshiro2010}).  
}\end{example}

In general, a kei has good involutions besides the identity map.  However, the dihedral quandle whose cardinarity  is an odd integer has only the identity map as a good involution.  Good involutions of dihedral quandles are classified in \cite{KamadaOshiro2010}.  

\begin{example}{\rm 
Let $X = {\rm cong}(G)$ be the conjugation quandle of a group $G$, which is $G$ itself as a set and the quandle operation is the conjugation; $x^y := y^{-1}x y$.  Let $\rho: X \to X$ be the invertion; 
$\rho(x) = x^{-1}$.  Then $(X, \rho)$ is a symmetric quandle (cf. \cite{Kamada2007}).   
}\end{example}

\begin{example}\label{doublecover}{\rm 
For a quandle (or a rack)  $X$, 
let $\overline{X}$ be a copy of $X$ and let $D(X)$ be the disjoint union of $X$ and $\overline{X}$.  By $\overline{x}$, we mean the element of $\overline{X}$ corresponding to an element $x$ of $X$.  Define a binary operations on $D(X)$ as follows: 
$$ 
\begin{array}{cc}
x^y := x^y \in X \subset D(X), \quad   &  
x^{\overline{y}} : = x^{y^{-1}} \in X \subset D(X), \quad     \\
\overline{x}^y := \overline{x^y} \in \overline{X} \subset D(X), \quad   &  
\overline{x}^{\overline{y}} : = \overline{x^{y^{-1}}} \in \overline{X} \subset D(X),   \\    
\end{array}
$$
where $x, y \in X$. 
Let $\rho: D(X) \to D(X)$ be the involution interchanging $x$ and $\overline{x}$ $(x \in X)$.  
Then $(D(X), \rho)$ is a symmetric quandle (or a symmetric rack)(cf. \cite{Kamada2007}). 

}\end{example}

\section{The full knot quandles and knot symmetric quandles}

Let $K$ be a  (PL and locally flat or smooth) proper $n$-submanifold of an $(n+2)$-manifold $W$.   
Let $N(K)$ be a tubular neighborhood of $K$ and $E(K)
= {\rm cl}(W \setminus N(K))$ the exterior of $K$, and take a base point $\ast$ of $E(K)$.   

We call a pair $(D, a)$ a {\it tadpole} of $K$ if 
$D$ is an oriented meridian disk of $K$ and $a$ is an arc in the exterior $E(K)$ starting from a point of $\partial D$ and ending at the base point $\ast$.   
Let 
$\widetilde{Q}(W, K, \ast) := \{ [(D, a)] \}$ be the quandle consisting of all homotopy classes of 
tadpoles with the operation $[(D_1, a_1)]^{[(D_2, a_2)]} := [(D_1, a_1 \cdot a_2^{-1} \cdot \partial D_2 \cdot a_2)]$.  

\begin{figure}[h]
\begin{center}
\includegraphics[scale=0.35]{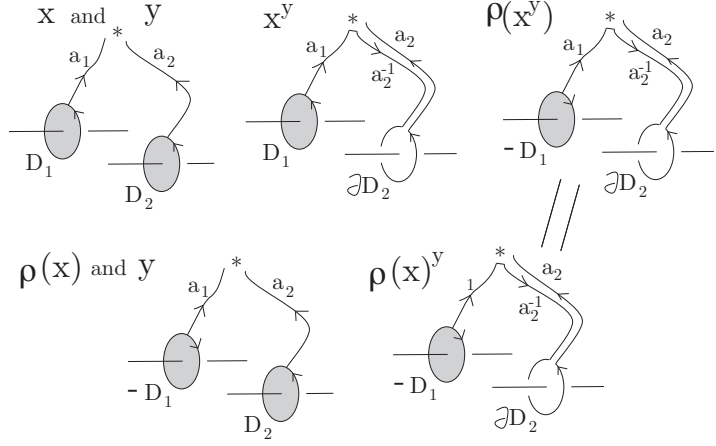}
\end{center}
\caption{$\rho(x^y) = \rho(x)^y$}
\label{fig01.eps}
\end{figure}

Let $\rho= \rho_K: \widetilde{Q}(W, K, \ast) \to \widetilde{Q}(W, K, \ast)$ be the map sending $[(D, a)]$ to $[(-D, a)]$, where $-D$ means $D$ with the reverse orientation.  It is easily verified that $\rho$ is a good involution (cf.  Figure~\ref{fig01.eps} for $\rho(x^y) = \rho(x)^y$), and we have a symmetric quandle $(\widetilde{Q}(W, K, \ast), \rho)$. 

When $W$ is connected, the isomorphism class of $\widetilde{Q}(W, K, \ast)$ does not depend on the base point $\ast$. 
So we denote it by $\widetilde{Q}(W, K)$.  

By an {\it $n$-manifold knot} in $\R^{n+2}$, we mean a closed (PL and locally flat or smooth) $n$-submanifold of $\R^{n+2}$.  

\begin{definition}[cf. \cite{Kamada2007, KamadaOshiro2010}]{\rm 
Let $K$ be a  (PL and locally flat or smooth) proper $n$-submanifold of an $(n+2)$-manifold $W$.   
The {\it full knot quandle} of $K$ is the quandle $\widetilde{Q}(W, K, \ast)$.  
The {\it knot symmetric quandle} (or the {\it fundamental symmetric quandle}) of $K$ is the symmetric quandle $(\widetilde{Q}(W, K, \ast), \rho)$, which we denote by $SQ(W, K, \ast)$.  When $K$ is an $n$-manifold knot in $\R^{n+2}$, 
the full knot quandle and the knot symmetric quandle of $K$ are denoted by 
$\widetilde{Q}(K)$  and $SQ(K)= (\widetilde{Q}(K), \rho)$, respectively.   
}\end{definition}

When $K$ is an orientable and oriented $n$-manifold knot in $\R^{n+2}$, if $D$ is an oriented meridian disk of $K$ whose orientation together with the orientation of $K$ matches the orientation of $\R^{n+2}$, we call a tadpole $(D, a)$ a {\it positive tadpole} of $K$.  The  {\it knot quandle} or the {\it positive knot quandle} of $K$ is the 
quandle $Q(K)$ consisting of all homotopy classes of positive tadpoles (cf. \cite{FennRourke1992, Joyce1982, Matveev1982}).  

\begin{remark}{\rm Let $K$ be an orientable and oriented $n$-manifold knot in $\R^{n+2}$.   
The knot quandle $Q(K)$ is a subquandle of the full knot quandle $\widetilde{Q}(K)$, and the knot symmetric quandle 
$SQ(K) = (\widetilde{Q}(K), \rho)$ is naturally identified with the double cover $D(Q(K))$ of the knot quandle $Q(K)$ in the sense of Example~\ref{doublecover}. 
}\end{remark}

\begin{figure}[h] 
\begin{center}
\includegraphics[scale=0.35]{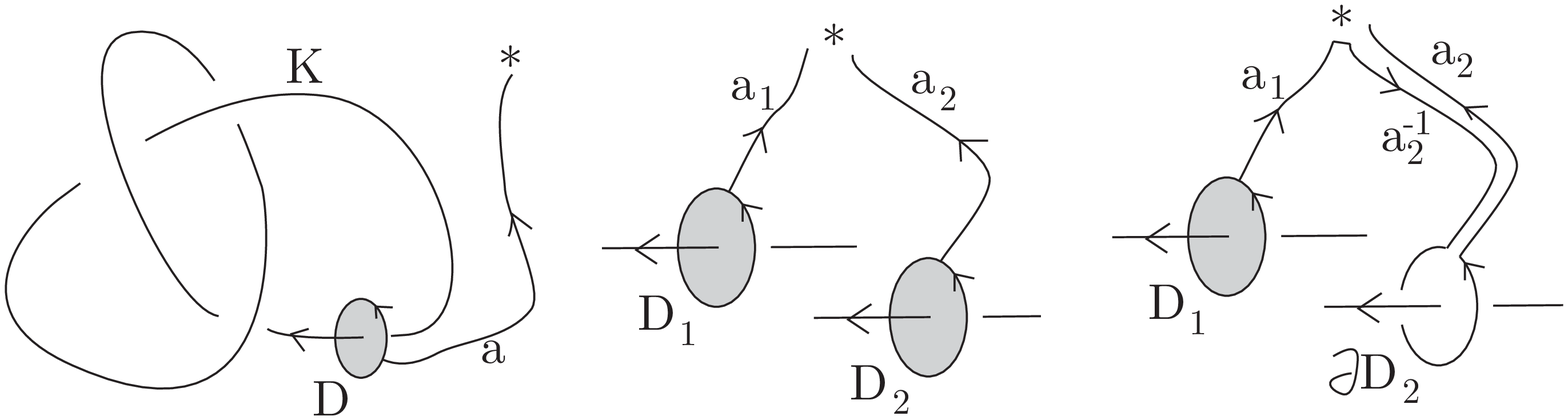}
\caption{The knot quandle and the operation}
\label{fig02.eps}
\end{center}
\end{figure}

\section{Presentations and associated groups} 

The notion of the associated group ${\rm As}(X)$ of a quandle or a rack $X$ is generalized to the notion of the 
associated group of a symmetric quandle or a symmetric rack.  

The {\it associated group} $ {\rm As}(X)$ of a quandle or a rack $X$ is defined as 
$$ {\rm As}(X)= \langle x  \in X \, | \,  x^y = y^{-1} x y  \, (x,y \in X) \rangle.$$

\begin{theorem}[\cite{FennRourke1992, Joyce1982, Matveev1982}]\label{thm:knotgrouppresentationA}
Let $K$ be an oriented $n$-manifold knot in  $\R^{n+2}$.  The associated group 
${\rm As}(Q(K))$ of the knot quandle $Q(K)$ is isomorphic to the knot group $G(K):= \pi_1(E(K))$.  
\end{theorem}

\begin{figure}[h]
\begin{center}
\includegraphics[scale=0.35]{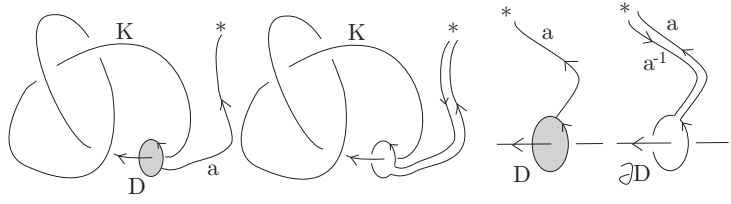}
\end{center}
\caption{The knot quandle $Q(K)$ to the knot group $G(K)$}
\label{fig03.eps}
\end{figure}

\begin{figure}[h]
\begin{center}
\includegraphics[scale=0.35]{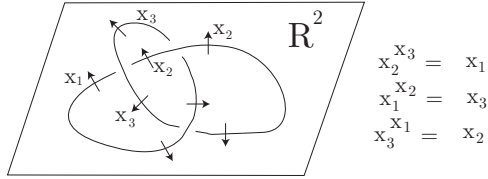}
\end{center}
\caption{A presentation from a diagram}
\label{fig04.eps}
\end{figure}

Let $K$ be an oriented knot in $\R^3$ and $D$ a diagram of $K$ in $\R^2$.  
As shown in \cite{FennRourke1992,Joyce1982,Matveev1982}, associated with the diagram $D$, we 
have a quandle $Q(D)$ defined by a presentation whose generators correspond to the arcs of the diagram and the relations correspond to the crossings.    
For example, 
$$Q(D) = \langle x_1, x_2, x_3 \,  |  \, x_2^{x_3} = x_1, x_1^{x_2} = x_3, x_3^{x_1} = x_2 \rangle_q$$
for the diagram in Figure~\ref{fig04.eps}.  
The notion of a presentation of a quandle is recalled in Section~\ref{sect:presentation}.  

\begin{theorem}[\cite{FennRourke1992,Joyce1982,Matveev1982}]\label{thm:knotqdlepresentationA}
Let $K$ be an oriented knot in $\R^3$ and $D$ a diagram of $K$ in $\R^2$.  
Then $Q(K)$ is isomorphic to $Q(D)$.  
\end{theorem}

Replacement of $x^y$ with $y^{-1}xy$ changes a presentation of a rack or a quandle 
$X$ to a group presentation of the associated group ${\rm As}(X)$ (cf. \cite{FennRourke1992}). 
For the previous example, we have 
$${\rm As}(Q(K))= G(K) = \langle x_1, x_2, x_3 \,  |  \, x_2^{x_3} = x_1, x_1^{x_2} = x_3, x_3^{x_1} = x_2 \rangle, $$
where $x^y$ stands for $y^{-1}xy$. This is a Wirtinger presentation of the knot group associated with the diagram. 

Theorem~\ref{thm:knotqdlepresentationA} is generalized to the higher dimensional case.  

Let $K$ be an $n$-manifold knot in $\R^{n+2}$.  We assume that it is in general position with respect to the projection 
$p: \R^{n+2} \to \R^{n+1}, \, (x_1, \cdots, x_{n+2}) \mapsto (x_1, \cdots, x_{n+1})$. The singular point set $\Delta$, that is the closure 
of the multiple point set $\{ y \in \R^{n+1} \, | \, 
\sharp( p^{-1}(y) \cap K ) > 1 \}$ in $\R^{n+1}$ (or its preimage in $K$), is  regarded as an $(n-1)$-dimensional stratified complex.  
It is well known that the $(n-1)$-dimensional strata consist of transverse double points (cf. \cite{RourkeSanderson1972}).  We call them the {\it double point strata} and denote by $\Delta^1$ the union of them.  
The lower dimensional strata are in general complicated.  
 (When $n=2$, the $0$-dimensional strata consists of triple points and branch points.)  For our purpose, classification of lower dimensional strata is not required at all since they do not contribute to the knot quandle. 
 
\begin{figure}[h]
\begin{center}
\includegraphics[scale=0.45]{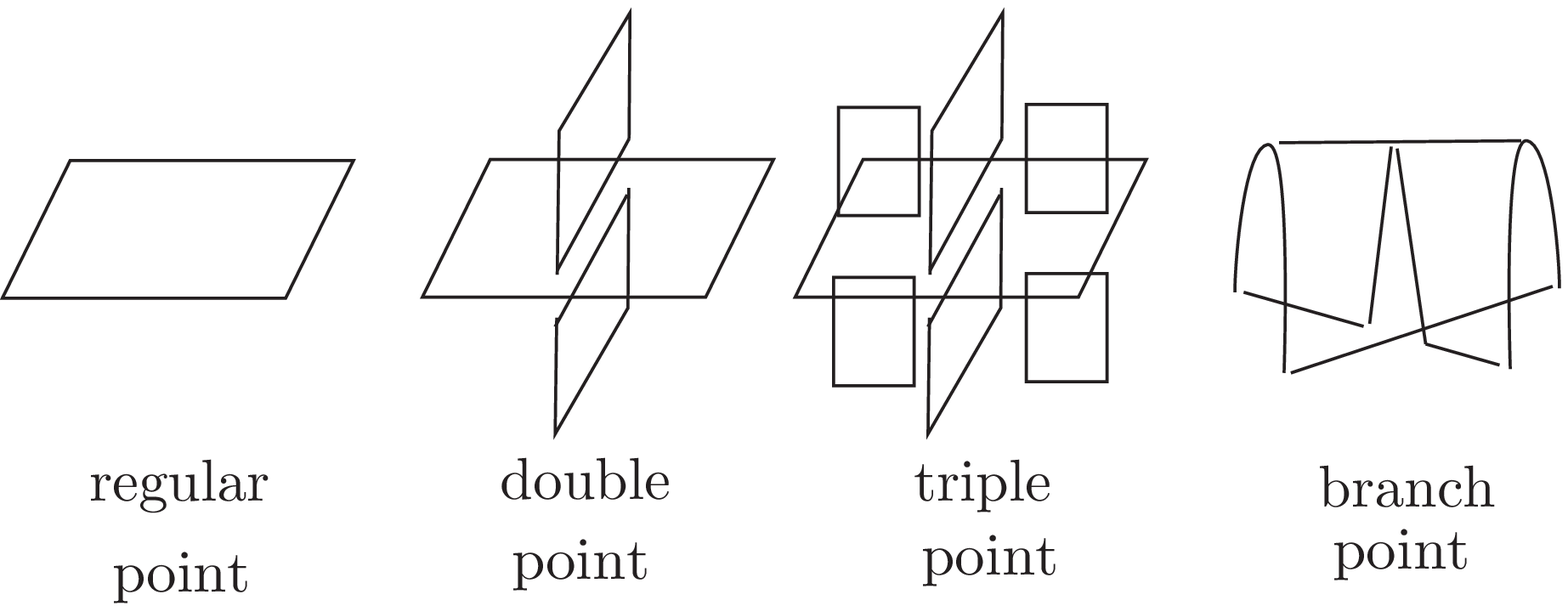}
\end{center}
\caption{Parts of a diagram of a surface knot}
\label{fig05.eps}
\end{figure}

A {\it diagram} of $K$ is a subset of $\R^{n+1}$ that is obtained from $p(K)$ by removing a regular neighborhood of the singular set $\Delta$ and put back the neighborhood of the upper branch of each double point stratum.  (See Figure~\ref{fig05.eps} for local pictures of a diagram when $n=2$.) 

The connected components of a diagram called the {\it sheets}.  (When $n=1$, they are also called the {\it arcs} of a diagram.) 

Now we assume that $K$ is orientable and oriented.  Then all sheets of the diagram $D$ are also oriented.  We usually present the orientation of each sheet by a normal vector in $\R^{n+1}$ as follows: Let $x$ be a point of a sheet and let $(v_1, \dots, v_n)$ be an $n$-tuple of tangent vectors at $x$ representing the orientation of the sheet.  Take a normal vector $w$ such that $(v_1, \dots, v_n, w)$ matches the orientation of $\R^{n+1}$.  Then we say that $w$ is the normal vector presenting the orientation of the sheet.

\begin{figure}[h]
\begin{center}
\includegraphics[scale=0.45]{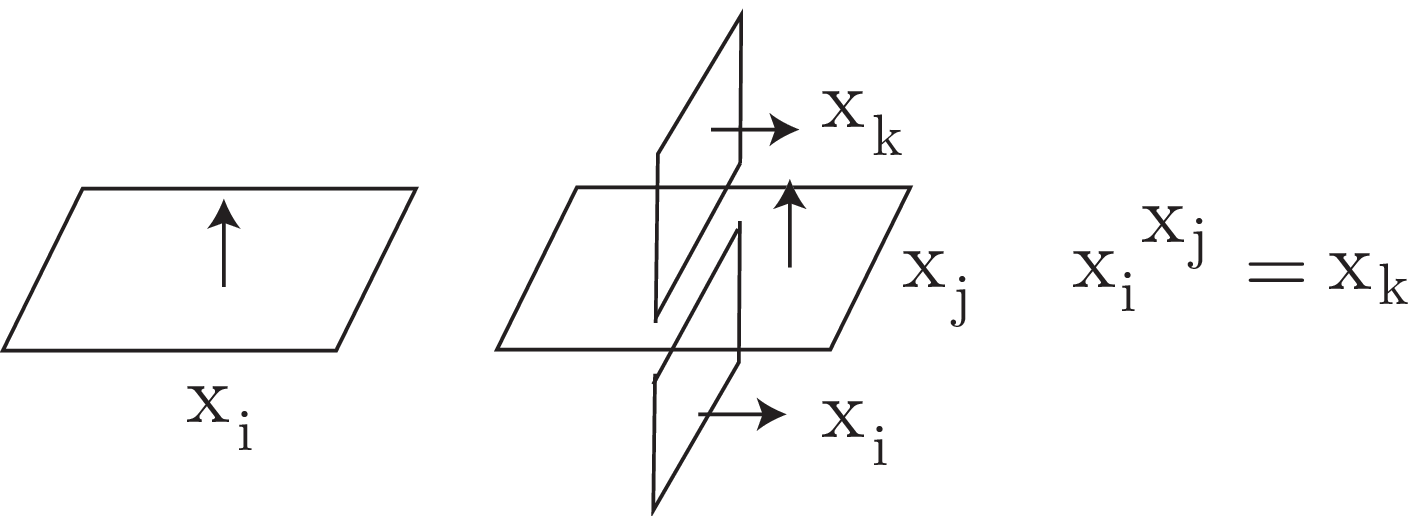}
\end{center}
\caption{Generators and relations}
\label{fig06.eps}
\end{figure}

Let $K$ be an oriented $n$-manifold knot in $\R^{n+2}$, and $D$ a diagram of $K$ in $\R^{n+1}$. 
Let $Q(D)$ be a quandle defined by a presentation whose generators correspond to the sheets of the diagram and the relations correspond to the double point strata. When the orientations of sheets are presented by normal vectors in $\R^{n+1}$ as in the left of Figure~\ref{fig06.eps},  the relation corresponding a double point stratum depicted in the right of Figure~\ref{fig06.eps} is $x_i^{x_j}=x_k$.  
(This definition is compatible with the definitions in \cite{CJKLS03, FennRourke1992, Joyce1982, Matveev1982} for the case where $K$ is an oriented knot in $\R^3$ or an oriented surface knot in $\R^4$.)    

\begin{theorem}\label{thm:knotqdlepresentationB}
Let $K$ be an oriented $n$-manifold knot in $\R^{n+2}$, and $D$ a diagram of $K$ in $\R^{n+1}$. Then 
$Q(K)$ is isomorphic to $Q(D)$.  
\end{theorem}

This theorem is proved in Section~\ref{sect:proof}. 

For a non-orientable $n$-manifold knot $K$ in $\R^{n+2}$, the knot quandle $Q(K)$ is not defined and 
we should consider the full knot quandle $\widetilde{Q}(K)$.  It will turn out that   
the knot symmetric quandle $SQ(K) = (\widetilde{Q}(K), \rho)$   is more natural and practical than the full knot quandle $\widetilde{Q}(K)$ itself.  

For a symmetric quandle or a symmetric rack, the associated group is defined as follows. 

\begin{definition}[\cite{Kamada2007, KamadaOshiro2010}]{\rm 
Let $(X, \rho)$ be a symmetric quandle or a symmetric rack.  
The {\it associated group} ${\rm As}(X, \rho)$ of $(X,\rho)$ is defined by 
$$ {\rm As}(X, \rho)= \langle x  \in X \, | \,  x^y = y^{-1} x y  \, (x,y \in X), \, 
\rho(x) = x^{-1} \, (x \in X) \rangle.$$
}\end{definition}

\begin{theorem}\label{thm:knotqdlepresentationGroup}
For an $n$-manifold knot $K$ of $\R^{n+2}$, the associated group ${\rm As}(\widetilde{Q}(K), \rho)$ of the knot symmetric quandle 
$(\widetilde{Q}(K), \rho)$ is isomorphic to the knot group 
$G(K)= \pi_1(\R^{n+2}\setminus K)$.  
\end{theorem}

We will explain how to get a presentation of the knot symmetric quandle 
$(\widetilde{Q}(K), \rho)$ from a diagram later.  
Theorem~\ref{thm:knotqdlepresentationGroup} is proved in Section~\ref{sect:proof}. 

\begin{example}{\rm 
Let $K$ be a standard 2-knot in $\R^4$.  The knot quandle and the full knot quandle are  
$$ 
{Q}(K) = \langle x \, | \,  \quad \rangle_q =\{x\} 
\quad \mbox{and} \quad 
\widetilde{Q}(K) = 
\langle x, y \, | \, x^y = x, y^x=y \rangle_q 
=\{x,y\}$$
where $x= [(D, a)]$ and $y =[(-D, a)]$.  They are a trivial quandle consisting of a single element and a trivial quandle consisting of $2$ elements, respectively. 
The associated groups of them are $\Z$ and $\Z \times \Z$.  On the other hand, the knot symmetric quandle   
$ SQ(K)$ is $(\widetilde{Q}(K), \rho)$ where 
the involution $\rho$ is given by $\rho(x)=y$. It has a presentation, as a symmetric quandle, 
$$ 
{SQ}(K) = \langle x \, | \, \quad \rangle_{sq},$$
and the associated group 
${\rm As}(\widetilde{Q}(K), \rho)$ has a group presentation 
$\langle x \rangle$, which is $\Z$. 
}\end{example}

\begin{figure}[h]
\begin{center}
\includegraphics[scale=0.45]{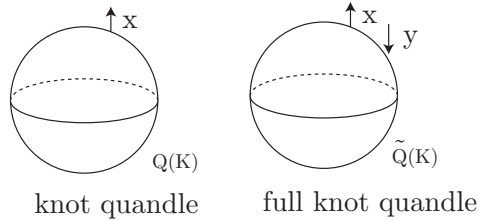}
\end{center}
\caption{The standard $2$-sphere}
\label{fig07.eps}
\end{figure}

As seen in this example and in Theorem~\ref{thm:knotgrouppresentationA}, for an oriented $n$-manifold knot, the knot quandle $Q(K)$ is sufficient to recover the knot group.  However this is not the case for a non-orientable $n$-manifold knot.  

\begin{example}{\rm 
Let $K$ be a standard projective plane in $\R^4$.  Since it is non-orientable, we cannot define the knot quandle $Q(K)$.  The full knot quandle is
$ \widetilde{Q}(K)=  \{ x \}$ and the knot symmetric quandle   
$ SQ(K)$ is $(\widetilde{Q}(K), \rho)$, where 
$\rho$ in the involution with $\rho(x)=x$.  
The knot symmetric quandle   
$ SQ(K)$ has a presentation, as a symmetric quandle, 
$$ 
{SQ}(K) = \langle x \, | \, \rho(x)= x \rangle_{sq},$$
and the associated group 
${\rm As}(\widetilde{Q}(K), \rho)$ has a group presentation 
$\langle x \, | \, x^{-1}= x \rangle$, which is $\Z/ 2\Z$. 
}\end{example}

Now we explain how to get a presentation of the knot symmetric quandle from a diagram.  

Let $K$ be an $n$-manifold knot in $\R^{n+2}$, and $D$ a diagram of $K$ in $\R^{n+1}$.  We assume that $K$ is in general position with respect to the projection 
$p: \R^{n+2} \to \R^{n+1}, \, (x_1, \cdots, x_{n+2}) \mapsto (x_1, \cdots, x_{n+1})$.  As stated before, a diagram is obtained from $p(K)$ by removing 
a regular neighborhood of the singular set $\Delta$ and by putting back the neighborhood of the upper branches of the double point strata.  
A {\it semi-sheet} of the diagram means a connected component of $p(K)$ removed the neighborhood of $\Delta$ (without putting back the neighborhood of the upper branches of the double point strata). 
Each semi-sheet is a compact {\it orientable} $n$-manifold in $\R^{n+1}$, even if $K$ is non-orientable (cf. Lemma~1 of \cite{Kamada2001}).  

\begin{figure}[h]
\begin{center}
\includegraphics[scale=0.45]{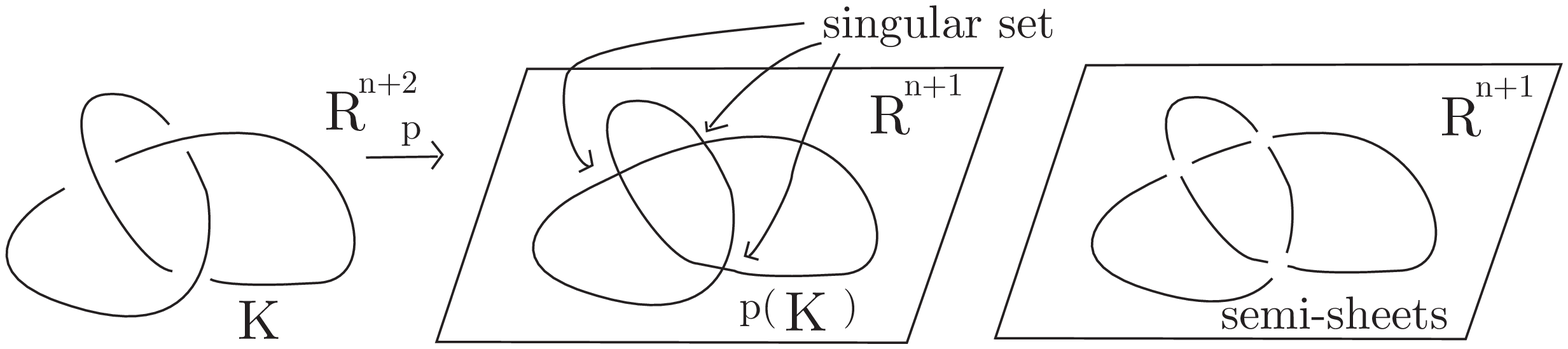}
\end{center}
\caption{The projection image $p(K) \subset \R^{n+1}$ and semi-sheets}
\label{fig08.eps}
\end{figure}

\begin{figure}[h]
\begin{center}
\includegraphics[scale=0.5]{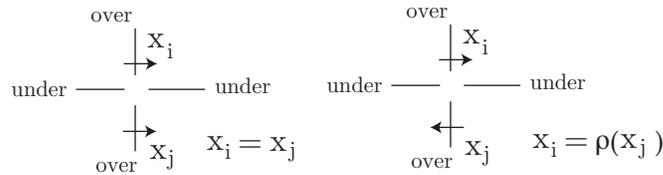}
\end{center}
\caption{$A$-relation}
\label{fig09.eps}
\end{figure}

\begin{figure}[h]
\begin{center}
\includegraphics[scale=0.5]{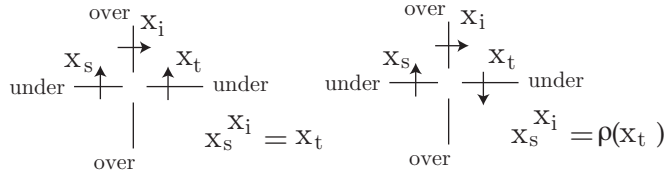}
\end{center}
\caption{$B$-relation}
\label{fig10.eps}
\end{figure}

Let $x_1, \cdots, x_m$ be the semi-sheets of $D$.  Give them normal vectors  
$w_1, \cdots, w_m$ arbitrarily.  We associate with each double point stratum an {\it $A$-relation} as in Fig.~\ref{fig09.eps} and a {\it $B$-relation} as in Fig.~\ref{fig10.eps} as follows: Let $x_i$ and $x_j$ be the upper semi-sheets.  If the normal vectors $w_i$ and $w_j$ of the upper semi-sheets $x_i$ and $x_j$ are  coherent (or incoherent, resp.), then let an $A$-relation be $x_i=x_j$ (or $x_i= \rho(x_j)$, resp.).  
Let $x_i$ be an upper semi-sheet, and let $x_s$ and $x_t$ be the lower semi-sheets such that the normal vector of $x_i$ is directed from $x_s$ to $x_t$.  If the normal vectors $w_s$ and $w_t$ of the lower semi-sheets $x_s$ and $x_t$ are   coherent (or incoherent, resp.), then 
let a $B$-relation be $x_s^{x_i}= x_t$ (or $x_s^{x_i}= \rho(x_t)$, resp.). 

\begin{definition}{\rm  
A {\it symmetric quandle presentation associated with the diagram $D$} is a presentation whose generators are $x_1, \cdots, x_m$ and the relations are 
$A$-relations and $B$-relations for double point strata.  We denote by $SQ(D)$ the symmetric quandle with this presentation. 
}\end{definition}

\begin{remark}{\rm   For a double point stratum where the normal vectors $w_i$ and $w_j$ of the upper semi-sheets are incoherent as depicted in the right of Figure~\ref{fig09.eps}, there are two possibilities for an $A$-relation: 
One is $x_i= \rho(x_j)$ and the other is $x_j= \rho(x_i)$.  Since $\rho$ is an involution, these two relations are equivalent.  

Similarly, there 
are two possibilities for a $B$-relation. However they are equivalent up to the $A$-relation.  
In the right of Figure~\ref{fig10.eps}, let $x_j$ be the 
over semi-sheet that is not labelled in the figure. 
For example, consider the case depicted in the right of the figure  
and assume that the normal vector of $x_j$ is directed from $x_t$ to $x_s$.  
Then we also have a relation $x_t^{x_j}=\rho(x_s)$ as 
a $B$-relation.    
The two relations $x_s^{x_i}=\rho(x_t)$ and  $x_t^{x_j}=\rho(x_s)$ are equivalent up to the $A$-relation;   
$(x_s^{x_i} = \rho(x_t))$  $\leftrightarrow$  $(x_s^{\rho(x_j)} = \rho(x_t))$  
 $\leftrightarrow$  $(x_s = \rho(x_t)^{x_j} )$
$\leftrightarrow$   $(\rho(x_s) = \rho( \rho(x_t)^{x_j} ))$
 $\leftrightarrow$   $(\rho(x_s) = x_t^{x_j})$.  

Thus it is sufficient for defining the symmetric quandle $SQ(D)$ to consider a single $A$-relation and a single $B$-relation for each double point stratum. 
}\end{remark}

\begin{example}{\rm 
Let $D$ be the diagram of a trefoil knot depicted in Figure~\ref{fig11.eps}. 
(The dotted arcs mean the neighborhood of the upper branches of the double point strata.) Choose normal vectors for the semi-sheets as in the figure. 
We may take $x_2=x_6$, $x_3=\rho(x_4)$ and $x_5=\rho(x_1)$ as $A$-relations, and  $x_3^{x_2}=x_1$, $x_2^{x_4}=\rho(x_5)$ and $x_4^{x_5}=\rho(x_6)$ as $B$-relations.  The symmetric quandle $SQ(D)$ associated with the diagram has a presentation 
\[
\langle 
x_1, \dots, x_6 \, | \, 
\begin{array}{ccc}
x_2=x_6,  &  x_3=\rho(x_4), &  x_5=\rho(x_1), \\
x_3^{x_2}=x_1,  & x_2^{x_4}=\rho(x_5),  &  x_4^{x_5}=\rho(x_6) 
\end{array}
\rangle_{sq}.
\]

\begin{figure}[h]
\begin{center}
\includegraphics[scale=0.45]{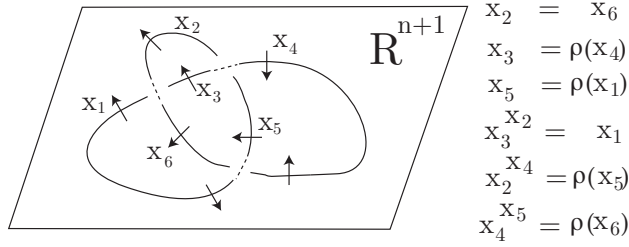}
\end{center}
\caption{Symmetric quandle presentation}
\label{fig11.eps}
\end{figure}

By eliminating $x_5$ and $x_6$ from this presentation, we have 
\[
\langle 
x_1, \dots, x_4 \, | \, 
\begin{array}{ccc}
    &  x_3=\rho(x_4), &   \\
x_3^{x_2}=x_1,  & x_2^{x_4}=x_1,  &  x_4^{\rho(x_1)}=\rho(x_2) 
\end{array}
\rangle_{sq}.
\]
The relation $x_3=\rho(x_4)$ is equivalent to $x_4=\rho(x_3)$, and the relation $x_2^{x_4}=x_1$ becomes $x_2^{\rho(x_3)}=x_1$, which is equivalent to $x_2=x_1^{x_3}$.  
The relation $x_4^{\rho(x_1)}=\rho(x_2)$ becomes $\rho(x_3)^{\rho(x_1)}=\rho(x_2)$, which is equivalent to $x_3^{\rho(x_1)}=x_2$ and hence to 
$x_3=x_2^{x_1}$.  
Thus we have a presentation 
\[
\langle 
x_1, \dots, x_3 \, | \, 
\begin{array}{ccc}
x_3^{x_2}=x_1,  & x_1^{x_3}=x_2,  &  x_2^{x_1}=x_3  
\end{array}
\rangle_{sq}.
\]
for $SQ(D)$. 
}\end{example}

\begin{theorem}\label{thm:knotqdlepresentationC}
Let $K$ be an  $n$-manifold knot in $\R^{n+2}$, and $D$ a diagram of $K$ in $\R^{n+1}$.  The knot symmetric quandle 
$SQ(K)$ of $K$ is isomorphic to the symmetric quandle $SQ(D)$ associated with the diagram $D$. 
\end{theorem}

We will prove this theorem in Section~\ref{sect:proof}. 

As a corollary, we see that the knot symmetric quandle $SQ(K)$ of an $n$-manifold knot in $\R^{n+2}$ has a finite presentation.  

Let $K$ be an $n$-manifold knot in $\R^{n+2}$. 
Let $(X, \rho)$ be a  symmetric quandle.  

\begin{definition}[\cite{Kamada2007, KamadaOshiro2010}]{\rm 
An {\it $(X, \rho)$-coloring} of a diagram of $K$ is an assignment of an element of $X$ to each semi-sheet of the diagram satisfying an A-relation and a B-relation  of each double point stratum.  
}\end{definition}

We denote by ${\rm Col}_{(X, \rho)}(D)$ the set of $(X, \rho)$-colorings of $D$, and by ${\rm Hom}(SQ(K), (X, \rho))$ the set of homomorphisms from the knot symmetric quandle $SQ(K)$ to the symmetric quandle $(X, \rho)$.

\begin{corollary}
There is a bijection from ${\rm Col}_{(X, \rho)}(D)$ to ${\rm Hom}(SQ(K), (X, \rho))$.  
In particular, the cardinarity $\sharp ({\rm Col}_{(X, \rho)}(D))$ is a knot invariant of $K$. 
\end{corollary}

Using a symmetric quandle, some invariants derived from a quandle, as quandle colorings, quandle cocycle invariants \cite{CJKLS03}, etc. can be also applied for non-orientable surface knots in $\R^4$ (cf. \cite{KamadaOshiro2010, Oshiro2010, Oshiro2011}).  Some other related researches are found in  \cite{CarterOshiroSaito2010, JangOshiro2012}.

\section{Presentations of symmetric quandles}
\label{sect:presentation}

This section is devoted to introducing the notion of a presentation of a symmetric quandle. 
First we recall the notions of a presentation of a rack and a presentation of a quandle due to Fenn and Rourke \cite{FennRourke1992}. 

Let $S$ be a non-empty set and let $F(S)$ be the free group generated by the elements of $S$. 
For an element $w \in F(S)$, we often denote $w^{-1}$ by $\overline{w}$.  

The {\it free rack} is a rack ${\rm FR}(S) := S \times F(S)$ with a rack operation 
$(a, w)^{(b, z)} = (a, w z^{-1} b z)$.  An element $(a,w)$ of ${\rm FR}(S)$ is also denoted by $a^w$; 
${\rm FR}(S) = \{ a^w \, | \, a \in S, w \in F(S) \}$.   Moreover we assume that the free group $F(S)$ acts on ${\rm FR(S)}$ from the right 
by ${\rm FR}(S) \times F(S) \to {\rm FR}(S); (a^w, z) \mapsto a^{wz}$.  Then ${\rm FR}(S)$ is regarded as an augmented rack with an augmentation $\partial: {\rm FR}(S) \to F(S)$  given by $a^w \mapsto w^{-1}aw$. (See page 355 of \cite{FennRourke1992} for the definition of an augmented rack.)   

Let $R$ be a subset of ${\rm FR}(S) \times {\rm FR}(S)$.  Let $\langle\langle R \rangle\rangle_r$ be 
the smallest congruence containing $R$ with respect to axioms of a rack.  

Start with the set $R$ and enlarge $R$ by repeating the following moves: 
\begin{itemize}
\item[(E1)] Add every diagonal element  $(x,x)$. 
\item[(E2)] If $(x,y) \in R$ then add $(y,x)$. 
\item[(E3)] If $(x,y) \in R$ and $(y,z)\in R$ then add $(x,z)$. 
\item[(R1)] If $(x,y) \in R$ then add $(x^a, y^a)$ and $(x^{\overline{a}}, y^{\overline{a}})$ 
for every $a \in S$.  
\item[(R2)] If $(x,y) \in R$ then add $(t^x, t^y)$ and $(t^{\overline{x}}, t^{\overline{y}})$ for every $t \in {\rm FR}(S)$.  
\end{itemize}
The result is $\langle\langle R \rangle\rangle_r$, which we call the {\it set of rack consequences} of $R$.  
Let 
$$ \langle S \, | \,  R \rangle_r := \frac{ {\rm FR}(S) }{ \langle\langle R \rangle\rangle_r}. $$
This is the rack whose generating set is $S$ and the defining relation set is $R$.  

For a rack $X$, if there is a map $f: S \to X$ which induces an isomorphism 
$ \langle S \, | \,  R \rangle_r \to X$ then we call it a {\it rack presentation}. 

For a quandle presentation, let $\langle\langle R \rangle\rangle_q$ be 
the smallest congruence containing $R$ with respect to axioms of a quandle.  

Start with the given set $R \subset {\rm FR}(S) \times {\rm FR}(S)$ and enlarge it by repeating the moves 
(E1), (E2), (E3), (R1), (R2) above and the following (Q). 
\begin{itemize}
\item[(Q)] Add $(a^a, a)$ for every $a \in S$.  
\end{itemize}
The result is  $\langle\langle R \rangle\rangle_q$, which we call the {\it set of quandle  consequences} of $R$.  
Let 
$$  \langle S \, | \,  R \rangle_q := \frac{ {\rm FR}(S) }{ \langle\langle R \rangle\rangle_q}. $$
This is the quandle whose generating set is $S$ and the defining relation set is $R$.  

For a quandle $X$, if there is a map $f: S \to X$ which induces an isomorphism 
$ \langle S \, | \,  R \rangle_q \to X$ then we call it a {\it quandle presentation}.

\begin{remark}{\rm 
Note that 
$\langle\langle R \rangle\rangle_q = \langle\langle R \cup 
\{ (a^a, a) \, | \, a \in S\} \rangle\rangle_r$, and 
$ \langle S \, | \,  R \rangle_q = \langle S \, | \,  R \cup \{ (a^a, a) \, | \, a \in S\}  \rangle_r$ (cf. Remark of p.~365 of \cite{FennRourke1992}).  
The {\it free quandle} on $S$ is $ \langle S \, | \,  \quad  \rangle_q$. 
}\end{remark}

Now we define presentations of a symmetric rack and a symmetric quandle.
First we define the free symmetric rack.  

Let $S$ be a non-empty set, and let $\overline{S}=\{ \overline{a} \, | \, a \in S\}$ be a copy of $S$ in which $\overline{a}$ stands for the copy of an element $a \in S$.  (We assume $S \cap \overline{S} = \emptyset$.)   

Let $F(S \cup \overline{S})$ be the free group generated by $S \cup \overline{S}$ and 
$\langle\langle a \overline{a} \, | \, a \in S \rangle\rangle$ the normal closure of 
$\{  a \overline{a} \, | \, a \in S \}$ in the free group.  The quotient 
$F(S \cup \overline{S}) / \langle\langle a \overline{a} \, | \, a \in S \rangle\rangle$ 
is naturally identified with the free group $F(S)$ by $a \mapsto a$ and $\overline{a} \mapsto a^{-1}$.  For an element $w \in F(S)$, we often denote $w^{-1}$ by $\overline{w}$.  

\begin{definition}{\rm 
The {\it free symmetric rack} on $S$ is a symmetric rack whose underlying rack is 
$$ {\rm FSR}(S) := S \times F(S) \amalg \overline{S} \times F(S) =  
(S \amalg \overline{S}) \times F(S) $$
with a rack operation 
$$ (a, w)^{(b,z)} := (a, w z^{-1} b z) $$
for $a, b \in S \amalg \overline{S}$ and $w, z \in F(S)$, 
and the involution $\rho$ is given by $\rho(a, w ) = (\overline{a}, w)$.   
We denote the free symmetric rack $({\rm FSR}(S), \rho)$ by ${\rm FSR}(S)$ for short. 
}\end{definition}

An element $(a, w) \in {\rm FSR}(S) $ is also denoted by $a^w$.  The image 
$\rho(a,w) = (\overline{a}, w)$ is denoted by $\rho(a^w)$ or $\overline{a^w}$, which is also 
by $\rho(a)^w$ or $\overline{a}^w$.  

The free group $F(S)$ acts on ${\rm FSR}(S)$ by ${\rm FSR}(S) \times F(S) \to {\rm FSR}(S); (a^w, z) \mapsto a^{wz}$.  Then ${\rm FSR}(S)$ is regarded as an augmented rack with an augmentation $\partial: {\rm FSR}(S) \to F(S)$  given by $a^w \mapsto w^{-1}aw$.  

\begin{remark}{\rm 
We define an {\it augmented symmetric rack} to be a symmetric rack $(X, \rho_X)$ whose underlying rack $X$ is an augmented rack with an augmentation $\partial: X \to G$ such that the involution $\rho_X$ is compatible with the canonical involution $\rho_G: G\to G; g \mapsto g^{-1}$ of the group $G$.  The free symmetric rack ${\rm FSR}(S) = ({\rm FSR}(S), \rho)$ with the action of $F(S)$ is an augmented symmetric rack.  
}\end{remark}
 
 Let $R$ be a subset of ${\rm FSR}(S) \times {\rm FSR}(S)$.  
 Let $\langle\langle R \rangle\rangle_{sr}$  be the smallest congruence containing $R$ with respect to axioms of a symmetric rack.  

Start with the set $R$ and enlarge $R$ by repeating the following moves: 
\begin{itemize}
\item[(SE1)] Add every diagonal element  $(x,x)$. 
\item[(SE2)] If $(x,y) \in R$ then add $(y,x)$. 
\item[(SE3)] If $(x,y) \in R$ and $(y,z)\in R$ then add $(x,z)$. 
\item[(SR1)] If $(x,y) \in R$ then add $(x^a, y^a)$ and $(x^{\overline{a}}, y^{\overline{a}})$ 
for every $a \in S$.  
\item[(SR2)] If $(x,y) \in R$ then add $(t^x, t^y)$ and $(t^{\overline{x}}, t^{\overline{y}})$ for every $t \in {\rm FSR}(S)$.  
\item[(S)] If $(x,y) \in R$ then add $(\overline{x}, \overline{y})$.  
\end{itemize}
The result is  $\langle\langle R \rangle\rangle_{sr}$, which we call the {\it set of symmetric rack consequences} of $R$.  
Let 
$$ \langle S \, | \, R \rangle_{sr} := \frac{ {\rm FSR}(S) }{ \langle\langle R \rangle\rangle_{sr}}. $$
This is the symmetric rack whose generating set is $S$ and the defining relation set is $R$.  

For a symmetric rack $X=(X, \rho)$, if there is a map $f: S \to X$ which induces an isomorphism 
$ \langle S \, | \, R \rangle_{sr} \to X$ then we call it a {\it symmetric rack presentation}. 

For a symmetric quandle presentation, let $\langle\langle R \rangle\rangle_{sq}$  be the smallest congruence containing $R$ with respect to axioms of a symmetric quandle.  

Start with the given set $R \subset {\rm FSR}(S) \times {\rm FSR}(S)$ and enlarge it by repeating the moves 
(SE1), (SE2), (SE3), (SR1), (SR2), (S) above and the move (Q): 
\begin{itemize}
\item[(Q)] Add $(a^a, a)$ for every $a \in S$.  
\end{itemize}
The result is $\langle\langle R \rangle\rangle_{sq}$, which we call the {\it set of symmetric quandle  consequences} of $R$.  
Let 
$$ \langle S \, | \, R \rangle_{sq} := \frac{ {\rm FSR}(S) }{ \langle\langle R \rangle\rangle_{sq}}. $$
This is the symmetric quandle whose generating set is $S$ and the defining relation set is $R$.

\begin{remark}{\rm 
Note that 
$\langle\langle R \rangle\rangle_{sq} = \langle\langle R \cup 
\{ (a^a, a) \, | \, a \in S\} \rangle\rangle_{sr}$ and that 
$ \langle S \, | \, R \rangle_{sq} = \langle S \, | \,  R \cup \{ (a^a, a) \, | \, a \in S\} \rangle_{sr}$.
}\end{remark}

\begin{definition}{\rm 
The {\it free symmetric quandle} on $S$ is $ \langle S \, | \, \quad \rangle_{sq}$.  We denote it by ${\rm FSQ}(S)$. 
}\end{definition}

\section{Proofs of theorems} 
\label{sect:proof}

In this section we prove 
Theorems $\ref{thm:knotqdlepresentationB}$,  $\ref{thm:knotqdlepresentationGroup}$ and 
$\ref{thm:knotqdlepresentationC}$.  

{\it Proof of Theorem~\ref{thm:knotqdlepresentationC}}.   
We use a similar argument as in the proof of Theorem~2 of \cite{Kamada2001}.  
Let $x_1, \dots, x_m$ be the semi-sheets of the diagram $D$, and let $w_1, \dots, w_m$ be normal vectors assigned for the semi-sheets.  
Let $N(\Delta)$ be a regular neighborhood of the singular set $\Delta$ and put 
$\Sigma := {\rm cl}(p(K) \setminus N(\Delta)) \subset \R^{n+1}$, which is the union of the semi-sheets.  

Without loss of generality, we may assume that $K$ in $\R^{n+2}= \R^{n+1} \times \R$ is in a form such that 
$$
K \cap \R^{n+1} \times [0, \infty)  =  (\partial \Sigma) \times [0, 1) \cup \Sigma \times \{ 1 \}  \\ 
$$
and the base point $\ast$ is a point in $\R^{n+2}= \R^{n+1} \times \R$ whose last coordinate is sufficiently large.

For each $i \in \{1, \dots, m\}$, consider an oriented small interval $a_i$ in $\R^{n+1}$ intersecting the interior of the  semi-sheet $x_i$ in the same direction with $w_i$, and let  $D_i = a_i \times [0,2] \subset \R^{n+1} \times \R = \R^{n+2}$, which is a meridian disk of $K$.  We give an orientation to $\partial D_i$ such that $\partial D_i \cap \R^{n+1} \times \{0\}$ is oriented in the same direction with $a_i$.  By this orientation, we assume that $D_i$ is an oriented meridian disk of $K$.  Let $b_i$ be a straight line in $\R^{n+2}$ starting from a point of $a_i \times \{2\} \subset \partial D_i$ to the base point $\ast$.  
We denote by the same symbol $x_i$ the element of $\widetilde{Q}(K)$ represented by the pair $(D_i, b_i)$, and also by the same symbol $x_i$ the element of the knot group $G(K)= \pi_1(\R^{n+2} \setminus K, \ast)$ represented by the loop $b_i^{-1} \cdot \partial D_i \cdot b_i$.   

Let $X=\{x_1, \dots, x_m\}$, and let ${\rm FSQ}(X)$  the free symmetric quandle on $X$ and $F(X)$ the free group on $X$. 
We can obtain a symmetric quandle homomorphism 
$$f_0: {\rm FSQ}(X) \to SQ(K)$$ with  
$x_i \mapsto x_i =[(D_i, b_i)]$ and $\overline{x_i} \mapsto \rho(x_i)=[(-D_i, b_i)]$ for $i \in \{1, \dots, m\}$, and a group 
homomorphism 
$$g_0: F(X) \to G(K)$$ with  
$x_i \mapsto x_i =[b_i^{-1} \cdot \partial D_i \cdot b_i]$ and $x_i^{-1} \mapsto x_i^{-1}=[b_i^{-1} \cdot \partial (-D_i) \cdot b_i]$ for $i \in \{1, \dots, m\}$.  

\begin{remark}{\rm 
Put $\R^{n+2}_+ = \R^{n+1} \times [0, \infty)$ and 
$K_+= K \cap \R^{n+2}_+$.  Then $K_+$ is a proper submanifold of $\R^{n+2}_+$.  
By the argument of Lemma~4 (2)  of \cite{Kamada2001}, we see that the knot group $\pi_1(\R^{n+2}_+ \setminus K_+, \ast)$ is identified with 
the free group $F(X)$ on $X$ and by using also the argument in \cite{Joyce1982}, we see that the 
knot symmetric quandle $SQ(\R^{n+2}_+, K_+, \ast)$ is identified with the free symmetric quandle ${\rm FSQ}(X)$ on $X$.  
}\end{remark}

We continue the proof.  
Put $\Delta= \Delta^1 \cup \Delta^2$, where $\Delta^1$ is the union of the double point strata and $\Delta^2$ is the union of the lower dimensional strata. 
We divide the regular neighborhood $N(\Delta)$ as follows: Let $N(\Delta^2)$ be a regular neighborhood of $\Delta^2$ in $\R^{n+1}$, and put $W_2= {\rm cl}(\R^{n+1} \setminus N(\Delta^2))$.  Let $N(\Delta_1)$ be a regular neighborhood of $\Delta^1 \cap W_2$, and put $W_1 = {\rm cl}(\R^{n+1} \setminus (N(\Delta^2) \cup N(\Delta^1))$.  We assume $N(\Delta)= N(\Delta^2) \cup N(\Delta^1)$.  

\begin{figure}[h]
\begin{center}
\includegraphics[scale=0.25]{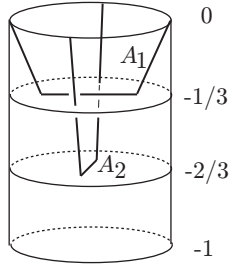}
\end{center}
\caption{Around the double point strata}
\label{fig12.eps}
\end{figure}

Let $A_1$ and $A_2$ be arcs in a cylinder $D^2 \times [-1, 0]$ as in Figure~\ref{fig12.eps}.  Note that $N(\Delta^1)$ is a trivial $D^2$-bundle over $\Delta^1 \cap W_2$, since each component of $\Delta^1$ has a trivialization determined from the four sheets around it.  Identify $N(\Delta^1)$ with $(\Delta^1 \cap W_2) \times D^2$ and $N(\Delta^1) \times [-1,0]$ with $(\Delta^1 \cap W_2) \times (D^2 \times [-1,0])$.  We may assume that $K$ restricted to $N(\Delta^1) \times [-1,0]$ is $(\Delta^1 \cap W_2) \times (A_1 \cup A_2)$ in 
$(\Delta^1 \cap W_2) \times (D^2 \times [-1,0])$.  Then we see that for each double point stratum, we have an $A$-relation and a $B$-relation.  
Therefore, the homomorphisms $f_0$ and $g_0$ induces 
a symmetric quandle homomorphism 
$$f : SQ(D) \to SQ(K)$$ with  
$x_i \mapsto x_i =[(D_i, b_i)]$ and $\overline{x_i} \mapsto \rho(x_i)=[(-D_i, b_i)]$ for $i \in \{1, \dots, m\}$, and a group 
homomorphism 
$$g: G(D) \to G(K)$$ with  
$x_i \mapsto x_i =[b_i^{-1} \cdot \partial D_i \cdot b_i]$ and $x_i^{-1} \mapsto x_i^{-1}=[b_i^{-1} \cdot \partial (-D_i) \cdot b_i]$ for $i \in \{1, \dots, m\}$.

In \cite{Kamada2001} it is proved by using the van Kampen theorem that the group homomorphism $g: G(D) \to G(K)$ is an isomorphism  from the fact that $\Delta^2$ has codimension more than two in $\R^{n+1}$.  A similar argument is applied for quandles (cf. \cite{Joyce1982}) and we see that the quandle homomorphism $f : SQ(D) \to SQ(K)$ is an isomorphism.  \qed 

\begin{remark}{\rm 
In the proof of Theorem~\ref{thm:knotqdlepresentationC}, we used the van Kampen theorem for quandles in \cite{Joyce1982} to obtain the fact that $f : SQ(D) \to SQ(K)$ is an isomorphism.  This fact can be also proved by a similar argument with the proof of Theorem~4.7 in \cite{FennRourke1992}.  
}\end{remark}

{\it Proof of Theorem~\ref{thm:knotqdlepresentationGroup}}.  
By Theorem~\ref{thm:knotqdlepresentationC}, the knot symmetric quandle $SQ(K)= (\widetilde{Q}(K), \rho)$ has a presentation associated with a diagram $D$.  In the proof of Theorem~\ref{thm:knotqdlepresentationC} (or by Theorem~2 of  \cite{Kamada2001}), we can obtain a group presentation of the knot group $G(K)$ associated with the diagram $D$, that is obtained from the symmetric quandle presentation associated with $D$ by replacing the quandle operation $x^y$ with conjugation $y^{-1}xy$, and $\rho(x)$ with $x^{-1}$.  This is the associated group of the symmetric quandle $SQ(K)= (\widetilde{Q}(K), \rho)$.  \qed 

{\it Proof of Theorem~\ref{thm:knotqdlepresentationB}}.   
In the proof of Theorem~\ref{thm:knotqdlepresentationC}, take the normal vectors $w_1, \dots, w_m$ such that they present the orientation of the semi-sheets.  Then all $A$-relations are in the case of $x_i=x_j$ as in the left of Figure~\ref{fig09.eps}.  Thus we can reduce the generating set from the set of semi-sheets to the set of sheets of the diagram.  The $B$-relations are as in the left of Figure~\ref{fig10.eps}. Then we see that the knot quandle $Q(K)$ has a quandle presentation associated with the oriented diagram $D$.  \qed 

\begin{remark}{\rm 
J{\'o}zef H. Przytycki and Witold Rosicki are studying quandles and 
cocycle invariants of $n$-manifold knots in $\R^{n+2}$ in \cite{PrzytyckiWitold2014}.   

}\end{remark}


\begin{thebibliography}{9} 


\bibitem{CJKLS03} J. S. Carter, D. Jelsovsky, S. Kamada, L. Langford and M. Saito, {\it Quandle cohomology and state-sum invariants of knotted curves and surfaces}, Trans. Amer. Math. Soc. {\bf 355} (2003), 3947--3989.

\bibitem{CarterOshiroSaito2010} J. S. Carter, K. Oshiro and M. Saito, {\it Symmetric extensions of dihedral quandles and triple points of non-orientable surfaces}, Topology Appl. {\bf 157} (2010), 857--869. 

\bibitem{Dehornoy2000} P. Dehornoy, {\it Braids and self-distributivity}, Progress in Mathematics Vol. 192, Birkhauser Verlag, Basel - Boston - Berlin, 2000. 

\bibitem{FennRourke1992} R. Fenn and C. Rourke, {\it Racks and links in codimension two}, J. Knot Theory Ramifications {\bf 1} (1992), 343--406.  

\bibitem{FRS2003} R. Fenn, C. Rourke and B. Sanderson, {\it The rack space}, Trans. Amer. Math. Soc. {\bf 359} (2007), 701--740.

\bibitem{JangOshiro2012} Y. Jang and K. Oshiro, {\it Symmetric quandle colorings for spatial graphs and handlebody-links}, J. Knot Theory Ramifications {\bf 21} (2012), 12500050, 16pp. 

\bibitem{Joyce1982} D. Joyce, {\it A classifying invariants of knots, the knot quandle}, J. Pure Appl. Alg. {\bf 23} (1982), 37--65.

\bibitem{Kamada2001} S. Kamada, {\it Wirtinger presentations for higher dimensional manifold knots obtained from diagrams}, Fund. Math. {\bf 168} (2001), 105--112.

\bibitem{Kamada2007} S. Kamada, {\it Quandles with good involutions, their homologies and knot invariants}, in: Intelligence of Low Dimensional Topology 2006, pp. 101--108, 
Ser. Knots Everything, Vol. 40, World Sci. Publ., Hackensack, NJ, 2007. 

\bibitem{KamadaOshiro2010} S. Kamada and K. Oshiro, {\it Homology groups of symmetric quandles and cocycle invariants of links and surface-links},  Trans. Amer. Math. Soc. {\bf 362} (2010), 5501--5527. 


\bibitem{Matveev1982} S. Matveev, {\it Distributive groupoids in knot theory}, 
Mat. Sb. (N.S.) {\bf 119 (161)} (1982) 78--88, 160; English translation: 
Math. USSR-Sb. {\bf 47} (1984), 73--83.


\bibitem{Oshiro2010} K. Oshiro, {\it Triple point numbers of surface-links and symmetric quandle cocycle invariants}, 
Algebr. Geom. Topol. {\bf 10} (2010), 853--865. 

\bibitem{Oshiro2011} K. Oshiro, {\it Homology groups of trivial quandles with good involutions and triple linking numbers of surface-links}, J. Knot Theory Ramifications {\bf 20} (2011), 595--608. 

\bibitem{PrzytyckiWitold2014} J. Przytycki and R. Witold, {\it Cocycle invariants of codimension $2$ embeddings of manifolds}, 
Banach Center Publications, to appear. 

\bibitem{RourkeSanderson1972} C. P. Rourke and B. J. Sanderson,   Introduction to Piecewise-Linear Topology, 
Ergeb. Math. Grenzgeb. 69, Springer, 1972. 

\bibitem{Takasaki1943} M. Takasaki, {\it Abstraction of symmetric transformations}, (in Japanese), Tohoku Math. J. {\bf 49} (1943), 145--207. 


\end{thebibliography}
\end{document}